\documentclass{llncs}
\usepackage{epsfig}
\newcommand{\R}{{I\!\!R}}
\newcommand{\N}{{I\!\!N}}

\def\R{{\rm I}\! {\rm R}}

\def\pt{\frac{\partial}{\partial t}}
\def\pvr{\frac{\partial}{\partial \vec r}}
\newcommand{\uu}[1]{\underline{\underline{#1}}}
\newcommand{\ttt}[1]{\mbox{\scriptsize #1}}
\renewcommand{\vec}[1]{\mbox{\boldmath $ #1 $}}

\newtheorem{algorithm}[theorem]{Algorithm}

\begin{document}

\pagestyle{headings}

\title{Numerical Methods of the Maxwell-Stefan Diffusion Equations and Applications in Plasma and Particle Transport}
\author{J\"urgen Geiser}
\institute{Ruhr University of Bochum, \\
The Institute of Theoretical Electrical Engineering, \\
Universit\"atsstra"se 150, D-44801 Bochum, Germany \\
\email{juergen.geiser@ruhr-uni-bochum.de}}
\maketitle

\begin{abstract}

In this paper, we present a model based on a local thermodynamic 
equilibrium, weakly ionized plasma-mixture model used for 
medical and technical applications in etching processes.
We consider a simplified model based on the Maxwell-Stefan model, 
which describe multicomponent diffusive fluxes in the gas mixture.
Based on additional conditions to the fluxes, we obtain an 
irreducible and quasi-positive diffusion matrix. Such problems 
results into nonlinear diffusion equations, which are more
delicate to solve as standard diffusion equations with 
Fickian's approach.
We propose explicit time-discretisation methods
embedded to iterative solvers for the nonlinearities.
Such a combination allows to solve the delicate nonlinear differential
equations more effective.
We present some first ternary component gaseous mixtures and discuss the
numerical methods.

\end{abstract}



\section{Introduction}

We are motivated to understand the gaseous mixtures of
a normal pressure and room temperature plasma.
The understanding of normal pressure, room temperature plasma applications
is important for applications in medical and technical processes.
Since many years, the increasing importance of plasma chemistry 
based on the multi-component plasma is a key factor to understand
the gaseous mixture processes, see for low pressure plasma \cite{sene06} 
and for atmospheric pressure regimes \cite{tanaka2004}.

We consider a simplified Maxwell-Stefan diffusion equation to model the
gaseous mixture of multicomponent Plasma.
While most classical description of the diffusion goes back to the Fickian's 
approach, see \cite{fick1855}, we apply the modern description of the
multicomponent diffusion based on the Maxwell-Stefan's approach, see \cite{maxwell1855}.
The novel approach considers are more detail description of the flux and 
concentration, which are in deed not only proportional coupled as in the
simplified Fickian's approach.
Here, we deal with a inter-species force balance, which allows to model cross-effects, e.g., so-called reverse diffusion (up-hill diffusion in direction of the gradients).

Such a more detailed modeling results in irreducible and quasi-positive 
diffusion-matrices, which can be reduced by transforming with 
reductions or transforming with Perron-Frobenius theorems to
solvable partial differential equations, see \cite{bothe2011}.
The obtained system of nonlinear partial differential equations are 
delicate to solve and numerically, we have taken into account
linearisation methods, e.g., iterative fix-point schemes.

The paper is outlined as follows.

In section \ref{modell} we present our mathematical model.
A possible reduced model for the further approximations
is derived in Section \ref{simply}.
In section \ref{method}, we discuss the underlying numerical schemes.
The first numerical results are presented in Section \ref{exper}.
In the contents, that are given in Section \ref{concl}, 
we summarize our results.

\section{Mathematical Model}
\label{modell}

For the full plasma model, we assumes that the neutral particles can be 
described as fluid dynamical model, where the 
elastic collision define the dynamics and few inelastic
collisions are, among other reasons, responsible for the 
chemical reactions.

To describe the individual mass densities, as
well as the global momentum and the global
energy as the dynamical conservation quantities
of the system, corresponding conservation equations
are derived from Boltzmann equations.

The individual character of each species 
is considered by mass-conservation equations and
the so-called difference equations.

The extension of the non-mixtured multicomponent transport
model, \cite{sene06} is done with respect to the
collision integrals related to the right-hind side sources of the
conservation laws.

The conservation laws of the neutral elements are given as
\begin{eqnarray*}
  \pt \rho_s + \pvr \cdot \rho_s \vec u_s &=& m_s Q_n^{(s)}, \\
  \pt \rho \vec u + \pvr \cdot \left(\uu{P}^* + \rho \vec u \vec u \right) &=& - Q_m^{(e)}, \\
  \pt \mathcal{E}_{\ttt{tot}}^*
     + \pvr \cdot \left(\mathcal{E}_{\ttt{tot}}^* \vec u + \vec q^* + \uu{P}^* \cdot \vec u \right)
     &=& - Q_{\mathcal{E}}^{(e)}.
\end{eqnarray*}
where $\rho_s$ : density of species $i$, $\rho = \sum_{i=1}^N \rho_i$, $\vec u$ : velocity,  $\mathcal{E}_{\ttt{tot}}^*$ : total energy of the neutral particles. \\

Further the variable $Q_n^{(s)}$ is the collision term of the mass 
conservation equation, $Q_m^{(e)}$ is the collision term of the 
momentum conservation equation and $ Q_{\mathcal{E}}^{(e)}$ is the collision 
term of the energy conservation equation.

We derive the collision term with respect to the
Chapmen-Enskog method, see \cite{chapman1990},  and achieve for the first derivatives the following 
results:
\begin{eqnarray}
&&  m_s Q_n^{(s)} = -  \nabla \cdot (\rho_i  \sum_{j=0} {\bf V}_i^{j} ), \\
&&  Q_m^{(e)} = - \sum_{i=1}^{n_s} \rho_i F_i \\
&& Q_{\mathcal{E}}^{(e)} = -  \sum_{i=1}^{n_s} \rho_i \rho F_i ({\bf u} + \sum_{j=0} V_i^{(j)}) ,
\end{eqnarray}
where $i = 1, \ldots, n_s$, $F_i$ is an external force per unit mass (see Boltzmann equation), further the diffusion velocity is given as:
\begin{eqnarray}
&& {\bf V}_i^{0} = 0 \\
&&  {\bf V}_i^{1} =  - \sum_{j=1}^N D_{ij} (d_j + k_{T_j} \frac{\Delta T}{T} ) ,
\end{eqnarray}
where $\sum_{i=1}^N d_i = 0$,
\begin{eqnarray}
d^∗_i = \nabla x_i + x_i \frac{\nabla p}{p} - \frac{\rho_i}{\rho} F_i , \\
d_i = d^∗_i - y_i \sum_{j} d_j^* ,
\end{eqnarray}
where $x_i = \frac{n_s}{n}$ is the molar fraction of species $i$.

We have an additional constraint based on the mass fraction of each species:
\begin{eqnarray}
&& \frac{\partial}{\partial t} y_i + \nabla y_i = R_i(y_1, \ldots, y_N) ,
\end{eqnarray}
where $y_i$ is the mass fraction of species $i$, $R_i$ is the net production
rate of species $i$ due to his reactions.

\begin{remark}
The full model problem consider a full coupled system of conservation
laws and Maxwell-Stefan equations. Each equations are coupled such that
the gaseous mixture influences the transport equations and 
vice verse. In the following, we decouple the equations system and
consider only the delicate Maxwell-Stefan equations.
\end{remark}

\section{Simplified Model with Maxwell-Stefan Diffusion Equations}
\label{simply}

We discuss in the following a multicomponent gaseous mixture with 
three species (ternary mixture). 
The model-problem is discussed in the experiments of 
Duncan and Toor, see \cite{duncan1962}.

Here, they studied an ideal gaseous mixture of the following components:
\begin{itemize}
\item Hydrogen ($H_2$, first species),
\item Nitrogen ($N_2$, second species),
\item Carbon dioxide ($CO_2$, third species).
\end{itemize}

The Maxwell-Stefan equations are given for the
three species as (see also \cite{boudin2012}):
\begin{eqnarray}
\label{part_1_eq_1}
&& \partial_t \xi_i + \nabla \cdot N_i = 0 , \; 1 \le i \le 3 , \\
\label{cond_1}
&& \sum_{j=1}^3 N_j = 0 , \\
\label{part_1_2}
&& \frac{\xi_2 N_1 - \xi_1 N_2}{D_{12}} + \frac{\xi_3 N_1 - \xi_1 N_3}{D_{13}} =  -  \nabla \xi_1 , \\
\label{part_1_3}
 && \frac{\xi_1 N_2 - \xi_2 N_1}{D_{12}} + \frac{\xi_3 N_2 - \xi_2 N_3}{D_{23}} =  -  \nabla \xi_2 ,
\end{eqnarray}
where the domain is given as $\Omega \in \R^d, d \in \N^+$ with $\xi_i \in C^2$.

For such ternary mixture, we can rewrite the three differential equations
(\ref{part_1_eq_1}) and (\ref{part_1_2})-(\ref{part_1_3}) with 
the help of the zero-condition (\ref{cond_1})
into two differential equations, given as:
\begin{eqnarray}
\label{part_2}
&& \partial_t \xi_i + \nabla \cdot N_i = 0 , \; 1 \le i \le 2 , \\
\label{part_2_1}
&& \frac{1}{D_{13}} N_1 + \alpha N_1 \xi_2 - \alpha N_2 \xi_1 =  -  \nabla \xi_1 , \\
\label{part_2_2}
 && \frac{1}{D_{23}} N_2 - \beta N_1 \xi_2 + \beta N_2 \xi_1 =  -  \nabla \xi_2 ,
\end{eqnarray}
where $\alpha = \left(\frac{1}{D_{12}} - \frac{1}{D_{13}}\right)$, 
$\beta = \left(\frac{1}{D_{12}} - \frac{1}{D_{23}}\right)$. 

Further we have the relations:
\begin{itemize}
\item Third mole-fraction: $\xi_3 = 1 - \xi_1 - \xi_2$, \\
\item Third molar flux: $N_3 = - N_1 - N_2$.
\end{itemize}

\section{Numerical Methods}
\label{method}

In the following, we discuss the numerical methods, which are 
based on iterative schemes with embedded explicit discretization 
schemes.
We apply the following methods:
\begin{itemize}
\item Iterative Scheme in time (Global Linearisation, Matrix Method),
\item Iterative Scheme in Time (Local Linearisation with Richardson's Method).
\end{itemize}

For the spatial discretization, we apply finite volume or finite difference
methods.
The underlying time-discretization is based on a first order explicit Euler
method.

\subsection{Iterative Scheme in time (Global Linearisation, Matrix Method)} 

We solve the iterative scheme:
\begin{eqnarray}
\label{ord_0}
&& \xi_{1}^{n+1} = \xi_1^n - \Delta t \; D_+ N_{1}^n , \\
&& \xi_{2}^{n+1} = \xi_2^n - \Delta t \; D_+ N_{2}^n , \\
&& \left( \begin{array}{c c}
A & B \\
C & D
\end{array} \right)
 \left( \begin{array}{l}
N_1^{n+1}  \\
N_2^{n+1} 
\end{array} \right) =
 \left( \begin{array}{l}
 - D_- \xi_1^{n+1} \\
 - D_- \xi_2^{n+1} 
\end{array} \right) 
\end{eqnarray}
for $j = 0, \ldots, J$ , where $\xi_1^n = (\xi_{1,0}^n, \ldots, \xi_{1, J}^n)^T$,
$\xi_2^n = (\xi_{2,0}^n, \ldots, \xi_{2, J}^n)^T$ and $I_J \in \R^{J+1} \times \R^{J+1}$,
 $N_1^n = (N_{1,0}^n, \ldots, N_{1, J}^n)^T$,
$N_2^n = (N_{2,0}^n, \ldots, N_{2, J}^n)^T$ and $I_J \in \R^{J+1} \times \R^{J+1}$,
where $n=0,1,2, \ldots, N_{end}$ and $N_{end}$ are the number of time-steps, i.d. $N_{end} = T / \Delta t$.

The matrices are given as:
\begin{eqnarray}
&& A, B, C, D \in \R^{J+1} \times \R^{J+1}, \\
&& A_{j,j} = \frac{1}{D_{13}} + \alpha \xi_{2,j} , \; j = 0 \ldots, J ,\\
&& B_{j,j} = - \alpha \xi_{1,j} , \; j = 0 \ldots, J , \\
&& C_{j,j} = - \beta \xi_{2, j} , \; j = 0 \ldots, J , \\
&& D_{j,j} =  \frac{1}{D_{23}}  + \beta \xi_{1,j}  , \; j = 0 \ldots, J ,\\
&& A_{i,j} = B_{i,j} =  C_{i,j} = D_{i,j} = 0 , \; i,j = 0 \ldots, J, \; i \neq J , 
\end{eqnarray}
means the diagonal entries given as for the scale case in
equation (\ref{part_2}) and the outer-diagonal entries are zero. \\
The explicit form with the time-discretization is given as:

\begin{algorithm}

1.) Initialisation $n=0$:

\begin{eqnarray}
\label{ord_0}
&& \left( \begin{array}{l}
N_1^{0}  \\
N_2^{0} 
\end{array} \right) =
 \left( \begin{array}{c c}
\tilde{A} & \tilde{B} \\
\tilde{C} & \tilde{D}
\end{array} \right)
 \left( \begin{array}{l}
 - D_- \xi_1^{0} \\
 - D_- \xi_2^{0} 
\end{array} \right) 
\end{eqnarray}
where $\xi_1^{0} = (\xi_{1,0}^{0}, \ldots, \xi_{1, J}^{0})^T$, $\xi_2^0 = (\xi_{2,0}^0, \ldots, \xi_{2, J}^0)^T$ and $\xi_{1,j}^{0} = \xi_1^{in}(j \Delta x), \; \xi_{2,j}^{0} = \xi_2^{in}(j \Delta x)$, $j = 0, \ldots, J$ and given as for the different intialisations, we have:
\begin{enumerate}
\item Uphill example
\begin{eqnarray}
\label{init}
&& \xi_1^{in}(x) = \left\{ \begin{array}{l l}
0.8 & \mbox{if} \; 0 \le x < 0.25 , \\
1.6 (0.75 - x) & \mbox{if} \; 0.25 \le x < 0.75 , \\
0.0 & \mbox{if} \; 0.75 \le x \le 1.0 , 
\end{array} \right. , \\
&& \xi_2^{in}(x) = 0.2 , \; \mbox{for all} \; x \in \Omega = [0,1] ,
\end{eqnarray}
\item Diffusion example (Asymptotic behavior)
\begin{eqnarray}
\label{init}
&& \xi_1^{in}(x) = \left\{ \begin{array}{l l}
0.8 & \mbox{if} \; 0 \le x  \in 0.5 , \\
0.0 & \mbox{else} , 
\end{array} \right. , \\
&&  \xi_2^{in}(x) = 0.2 ,  \; \mbox{for all} \; x \in \Omega = [0,1] ,
\end{eqnarray}
\end{enumerate}

The inverse matrices are given as:
\begin{eqnarray}
&& \tilde{A}, \tilde{B}, \tilde{C}, \tilde{D} \in \R^{J+1} \times \R^{J+1}, \\
&& \tilde{A}_{j,j} = \gamma_j (\frac{1}{D_{23}}  + \beta \xi_{1,j}^{0}) , \; j = 0 \ldots, J ,\\
&& B_{j,j} = \gamma_j \;  \alpha \xi_{1,j}^{0} , \; j = 0 \ldots, J , \\
&& C_{j,j} = \gamma_j \;  \beta \xi_{2, j}^{0} , \; j = 0 \ldots, J , \\
&& D_{j,j} = \gamma_j  (\frac{1}{D_{13}} + \alpha \xi_{2,j}^{0})  , \; j = 0 \ldots, J ,\\
&& \gamma_j = \frac{D_{13} D_{23}}{1 + \alpha D_{13} \xi_{2,j}^{0} + \beta D_{23} \xi_{1,j}^{0}} ,  \; j = 0 \ldots, J , \\
&& \tilde{A}_{i,j} = \tilde{B}_{i,j} =  \tilde{C}_{i,j} = \tilde{D}_{i,j} = 0 , \; i,j = 0 \ldots, J, \; i \neq J , 
\end{eqnarray}

Further the values of the first and the last grid points of $N$ are zero,
means $N_{1,0}^{0} = N_{1,J}^{0} = N_{2,0}^{0} = N_{2,J}^{0} = 0$ (boundary condition).

2.) Next time-steps (till $n = N_{end}$ ): \\

2.1) Computation of $\xi_1^{n+1}$ and $\xi_2^{n+1}$ 
\begin{eqnarray}
\label{ord_0}
&& \xi_{1}^{n+1} = \xi_1^n - \Delta t \; D_+ N_{1}^n , \\
&& \xi_{2}^{n+1} = \xi_2^n - \Delta t \; D_+ N_{2}^n ,
\end{eqnarray}

2.2) Computation of $N_1^{n+1}$ and $N_2^{n+1}$ 

\begin{eqnarray}
&& \left( \begin{array}{l}
N_1^{n+1}  \\
N_2^{n+1} 
\end{array} \right) =
 \left( \begin{array}{c c}
\tilde{A} & \tilde{B} \\
\tilde{C} & \tilde{D}
\end{array} \right)
 \left( \begin{array}{l}
 - D_- \xi_1^{n+1} \\
 - D_- \xi_2^{n+1} 
\end{array} \right) 
\end{eqnarray}
where $\xi_1^{n} = (\xi_{1,0}^{n}, \ldots, \xi_{1, J}^{n})^T$, $\xi_2^n = (\xi_{2,0}^n, \ldots, \xi_{2, J}^n)^T$ 
and the inverse matrices are given as:
\begin{eqnarray}
&& \tilde{A}, \tilde{B}, \tilde{C}, \tilde{D} \in \R^{J+1} \times \R^{J+1}, \\
&& \tilde{A}_{j,j} = \gamma_j (\frac{1}{D_{23}}  + \beta \xi_{1,j}^{n+1}) , \; j = 0 \ldots, J ,\\
&& B_{j,j} = \gamma_j \;  \alpha \xi_{1,j}^{n+1} , \; j = 0 \ldots, J , \\
&& C_{j,j} = \gamma_j \;  \beta \xi_{2, j}^{n+1} , \; j = 0 \ldots, J , \\
&& D_{j,j} = \gamma_j  (\frac{1}{D_{13}} + \alpha \xi_{2,j}^{n+1})  , \; j = 0 \ldots, J ,\\
&& \gamma_j = \frac{D_{13} D_{23}}{1 + \alpha D_{13} \xi_{2,j}^{n+1} + \beta D_{23} \xi_{1,j}^{n+1}} ,  \; j = 0 \ldots, J , \\
&& \tilde{A}_{i,j} = \tilde{B}_{i,j} =  \tilde{C}_{i,j} = \tilde{D}_{i,j} = 0 , \; i,j = 0 \ldots, J, \; i \neq J .
\end{eqnarray}

Further the values of the first and the last grid points of $N$ are zero,
means $N_{1,0}^{n} = N_{1,J}^{n} = N_{2,0}^{n} = N_{2,J}^{n} = 0$ (boundary condition).

3.) Do  $n = n+1$ and goto 2.) 

\end{algorithm}

\subsection{Iterative Scheme in Time (Local Linearisation with Richardson's Method}
We solve the iterative scheme given in the Richardson iterative scheme:
\begin{eqnarray}
\label{ord_0}
&& \xi_{1}^{n+1, k} = \xi_1^n - \Delta t \; D_+ N_{1}^{n+1} , \\
&& \xi_{2}^{n+1, k} = \xi_2^n - \Delta t \; D_+ N_{2}^{n+1} , \\
&& \left( \begin{array}{c c}
A^{n+1, k-1} & B^{n+1, k-1} \\
C^{n+1, k-1} & D^{n+1, k-1}
\end{array} \right)
 \left( \begin{array}{l}
N_1^{n+1}  \\
N_2^{n+1} 
\end{array} \right) =
 \left( \begin{array}{l}
 - D_- \xi_1^{n+1, k-1} \\
 - D_- \xi_2^{n+1, k-1} 
\end{array} \right) 
\end{eqnarray}
for $j = 0, \ldots, J$ , where $\xi_1^n = (\xi_{1,0}^n, \ldots, \xi_{1, J}^n)^T$,
$\xi_2^n = (\xi_{2,0}^n, \ldots, \xi_{2, J}^n)^T$ and $I_J \in \R^{J+1} \times \R^{J+1}$,
 $N_1^n = (N_{1,0}^n, \ldots, N_{1, J}^n)^T$,
$N_2^n = (N_{2,0}^n, \ldots, N_{2, J}^n)^T$ and $I_J \in \R^{J+1} \times \R^{J+1}$,
where $n=0,1,2, \ldots, N_{end}$ and $N_{end}$ are the number of time-steps, i.d. $N_{end} = T / \Delta t$.

Further $k = 1, 2, \ldots, K$ is the iteration index with \\
 where $\xi_1^{n+1, 0} = (\xi_{1,0}^n, \ldots, \xi_{1, J}^n)^T$,
$\xi_2^{n+1, 0} = (\xi_{2,0}^n, \ldots, \xi_{2, J}^n)^T$ and $I_J \in \R^{J+1} \times \R^{J+1}$ is the start solution given with the solution at $t = t^n$.

The matrices are given as:
\begin{eqnarray}
&& A^{n+1, k-1}, B^{n+1, k-1}, C^{n+1, k-1}, D^{n+1, k-1} \in \R^{J+1} \times \R^{J+1}, \\
&& A_{j,j}^{n+1, k-1} = \frac{1}{D_{13}} + \alpha \xi_{2,j}^{n+1, k-1} , \; j = 0 \ldots, J ,\\
&& B_{j,j}^{n+1, k-1} = - \alpha \xi_{1,j}^{n+1, k-1} , \; j = 0 \ldots, J , \\
&& C_{j,j}^{n+1, k-1} = - \beta \xi_{2, j}^{n+1, k-1} , \; j = 0 \ldots, J , \\
&& D_{j,j}^{n+1, k-1} =  \frac{1}{D_{23}}  + \beta \xi_{1,j}^{n+1, k-1}  , \; j = 0 \ldots, J ,\\
&& A_{i,j}^{n+1, i-1} = B_{i,j}^{n+1, i-1} =  C_{i,j}^{n+1, i-1} = D_{i,j}^{n+1, i-1} = 0 , \; i,j = 0 \ldots, J, \; i \neq J , 
\end{eqnarray}
means the diagonal entries given as for the scale case in
equation (\ref{part_2}) and the outer-diagonal entries are zero. \\
The explicit form with the time-discretization is given as:

\begin{algorithm}

1.) Initialisation $n=0$ with an explicit time-step (CFL condition is given):

\begin{eqnarray}
\label{ord_0}
&& \left( \begin{array}{l}
N_1^{0}  \\
N_2^{0} 
\end{array} \right) =
 \left( \begin{array}{c c}
\tilde{A} & \tilde{B} \\
\tilde{C} & \tilde{D}
\end{array} \right)
 \left( \begin{array}{l}
 - D_- \xi_1^{0} \\
 - D_- \xi_2^{0} 
\end{array} \right) 
\end{eqnarray}
where $\xi_1^{0} = (\xi_{1,0}^{0}, \ldots, \xi_{1, J}^{0})^T$, $\xi_2^0 = (\xi_{2,0}^0, \ldots, \xi_{2, J}^0)^T$ and $\xi_{1,j}^{0} = \xi_1^{in}(j \Delta x), \; \xi_{2,j}^{0} = \xi_2^{in}(j \Delta x)$, $j = 0, \ldots, J$ and given as for the different intialisations, we have:
\begin{enumerate}
\item Uphill example
\begin{eqnarray}
\label{init}
&& \xi_1^{in}(x) = \left\{ \begin{array}{l l}
0.8 & \mbox{if} \; 0 \le x < 0.25 , \\
1.6 (0.75 - x) & \mbox{if} \; 0.25 \le x < 0.75 , \\
0.0 & \mbox{if} \; 0.75 \le x \le 1.0 , 
\end{array} \right. , \\
&& \xi_2^{in}(x) = 0.2 , \; \mbox{for all} \; x \in \Omega = [0,1] ,
\end{eqnarray}
\item Diffusion example (Asymptotic behavior)
\begin{eqnarray}
\label{init}
&& \xi_1^{in}(x) = \left\{ \begin{array}{l l}
0.8 & \mbox{if} \; 0 \le x  \in 0.5 , \\
0.0 & \mbox{else} , 
\end{array} \right. , \\
&&  \xi_2^{in}(x) = 0.2 ,  \; \mbox{for all} \; x \in \Omega = [0,1] ,
\end{eqnarray}
\end{enumerate}

The inverse matrices are given as:
\begin{eqnarray}
&& \tilde{A}, \tilde{B}, \tilde{C}, \tilde{D} \in \R^{J+1} \times \R^{J+1}, \\
&& \tilde{A}_{j,j} = \gamma_j (\frac{1}{D_{23}}  + \beta \xi_{1,j}^{0}) , \; j = 0 \ldots, J ,\\
&& B_{j,j} = \gamma_j \;  \alpha \xi_{1,j}^{0} , \; j = 0 \ldots, J , \\
&& C_{j,j} = \gamma_j \;  \beta \xi_{2, j}^{0} , \; j = 0 \ldots, J , \\
&& D_{j,j} = \gamma_j  (\frac{1}{D_{13}} + \alpha \xi_{2,j}^{0})  , \; j = 0 \ldots, J ,\\
&& \gamma_j = \frac{D_{13} D_{23}}{1 + \alpha D_{13} \xi_{2,j}^{0} + \beta D_{23} \xi_{1,j}^{0}} ,  \; j = 0 \ldots, J , \\
&& \tilde{A}_{i,j} = \tilde{B}_{i,j} =  \tilde{C}_{i,j} = \tilde{D}_{i,j} = 0 , \; i,j = 0 \ldots, J, \; i \neq J , 
\end{eqnarray}

Further the values of the first and the last grid points of $N$ are zero,
means $N_{1,0}^{0} = N_{1,J}^{0} = N_{2,0}^{0} = N_{2,J}^{0} = 0$ (boundary condition).

2.) Next timesteps (till $n = N_{end}$ ) (iterative scheme restricted via the CFL condition based on the previous iterative solutions in the matrices): \\

2.1) Computation of $\xi_1^{n+1, I}$ and $\xi_2^{n+1, I}$ 
\begin{eqnarray}
\label{ord_0}
&& \xi_{1}^{n+1, k} = \xi_1^n - \Delta t \; D_+ N_{1}^{n+1} , \\
&& \xi_{2}^{n+1, k} = \xi_2^n - \Delta t \; D_+ N_{2}^{n+1} ,
\end{eqnarray}

2.2) Computation of $N_1^{n+1, k-1}$ and $N_2^{n+1, k-1}$ 

\begin{eqnarray}
&& \left( \begin{array}{l}
N_1^{n+1}  \\
N_2^{n+1} 
\end{array} \right) =
 \left( \begin{array}{c c}
\tilde{A}^{n+1, k-1} & \tilde{B}^{n+1, k-1} \\
\tilde{C}^{n+1, k-1} & \tilde{D}^{n+1, k-1}
\end{array} \right)
 \left( \begin{array}{l}
 - D_- \xi_1^{n+1, k-1} \\
 - D_- \xi_2^{n+1, k-1} 
\end{array} \right) 
\end{eqnarray}
where $\xi_1^{n} = (\xi_{1,0}^{n}, \ldots, \xi_{1, J}^{n})^T$, $\xi_2^n = (\xi_{2,0}^n, \ldots, \xi_{2, J}^n)^T$ 
and the inverse matrices are given as:
\begin{eqnarray}
&& \tilde{A}^{n+1, k-1}, \tilde{B}^{n+1, k-1}, \tilde{C}^{n+1, k-1}, \tilde{D}^{n+1, k-1} \in \R^{J+1} \times \R^{J+1}, \\
&& \tilde{A}_{j,j}^{n+1, k-1} = \gamma_j (\frac{1}{D_{23}}  + \beta \xi_{1,j}^{n+1, k-1}) , \; j = 0 \ldots, J ,\\
&& B_{j,j}^{n+1, k-1} = \gamma_j \;  \alpha \xi_{1,j}^{n+1, k-1} , \; j = 0 \ldots, J , \\
&& C_{j,j}^{n+1, k-1} = \gamma_j \;  \beta \xi_{2, j}^{n+1, k-1} , \; j = 0 \ldots, J , \\
&& D_{j,j}^{n+1, k-1} = \gamma_j  (\frac{1}{D_{13}} + \alpha \xi_{2,j}^{n+1, k-1})  , \; j = 0 \ldots, J ,\\
&& \gamma_j = \frac{D_{13} D_{23}}{1 + \alpha D_{13} \xi_{2,j}^{n+1, k-1} + \beta D_{23} \xi_{1,j}^{n+1, k-1}} ,  \; j = 0 \ldots, J , \\
&& \tilde{A}_{i,j}^{n+1, k-1} = \tilde{B}_{i,j}^{n+1, k-1} =  \tilde{C}_{i,j}^{n+1, k-1} = \tilde{D}_{i,j}^{n+1, k-1} = 0 , \; i,j = 0 \ldots, J, \; i \neq J .
\end{eqnarray}

Further the values of the first and the last grid points of $N$ are zero,
means $N_{1,0}^{n+1} = N_{1,J}^{n+1} = N_{2,0}^{n+1} = N_{2,J}^{n+1} = 0$ (boundary condition).

Further $k = 1, 2, \ldots, K$ is the iteration index with \\
 where $\xi_1^{n+1, 0} = (\xi_{1,0}^n, \ldots, \xi_{1, J}^n)^T$,
$\xi_2^{n+1, 0} = (\xi_{2,0}^n, \ldots, \xi_{2, J}^n)^T$ and $I_J \in \R^{J+1} \times \R^{J+1}$ is the start solution given with the solution at $t = t^n$.

3.) Do  $n = n+1$ and goto 2.) 

\end{algorithm}

\section{Numerical Experiments}
\label{exper}

In the following, we concentrate on the following three component system,
which is given as:
\begin{eqnarray}
\label{ord_0}
&& \partial_t \xi_i + \partial_x N_i = 0 , \; 1 \le i \le 3 , \\
&& \sum_{j=1}^3 N_j = 0 , \\
&& \frac{\xi_2 N_1 - \xi_1 N_2}{D_{12}} + \frac{\xi_3 N_1 - \xi_1 N_3}{D_{13}} =  -  \partial_x \xi_1 , \\
 && \frac{\xi_1 N_2 - \xi_2 N_1}{D_{12}} + \frac{\xi_3 N_2 - \xi_2 N_3}{D_{23}} =  -  \partial_x \xi_2 ,
\end{eqnarray}
where the domain is given as $\Omega \in \R^d, d \in \N^+$ with $\xi_i \in C^2$.

The parameters and the initial and boundary conditions are given as:
\begin{itemize}
\item $D_{12} = D_{13} = 0.833$ (means $\alpha = 0$) and $D_{23} = 0.168$ (Uphill diffusion, semi-degenerated Duncan and Toor experiment),
\item $D_{12} = 0.0833, D_{13} = 0.680$ and $D_{23} = 0.168$ (asymptotic behavior, Duncan and Toor experiment, see \cite{duncan1962}),
\item $J = 140$ (spatial grid points),
\item The time-step-restriction for the explicit method is given as: \\
 $\Delta t \le \frac{(\Delta x)^2}{2 \max\{D_{12}, D_{13}, D_{23}\}}$,
\item The spatial domain is $\Omega = [0, 1]$, the time-domain $[0, T] = [0, 1]$,
\item The initial conditions are:
\begin{enumerate}
\item Uphill example
\begin{eqnarray}
\label{init}
&& \xi_1^{in}(x) = \left\{ \begin{array}{l l}
0.8 & \mbox{if} \; 0 \le x < 0.25 , \\
1.6 (0.75 - x) & \mbox{if} \; 0.25 \le x < 0.75 , \\
0.0 & \mbox{if} \; 0.75 \le x \le 1.0 , 
\end{array} \right. , \\
&& \xi_2^{in}(x) = 0.2 , \; \mbox{for all} \; x \in \Omega = [0,1] ,
\end{eqnarray}
\item Diffusion example (Asymptotic behavior)
\begin{eqnarray}
\label{init}
&& \xi_1^{in}(x) = \left\{ \begin{array}{l l}
0.8 & \mbox{if} \; 0 \le x  \in 0.5 , \\
0.0 & \mbox{else} , 
\end{array} \right. , \\
&&  \xi_2^{in}(x) = 0.2 ,  \; \mbox{for all} \; x \in \Omega = [0,1] .
\end{eqnarray}
\end{enumerate}

\item The boundary conditions are of no-flux type:
\begin{eqnarray}
\label{init}
&& N_1 = N_2 = N_3 = 0 , \mbox{on} \; \partial \Omega \times [0,1] .
\end{eqnarray}
\end{itemize}

We could reduce to a simpler model problem as:
\begin{eqnarray}
\label{ord_0}
&& \partial_t \xi_i + \partial_x \cdot N_i = 0 , \; 1 \le i \le 2 , \\
&& \frac{1}{D_{13}} N_1 + \alpha N_1 \xi_2 - \alpha N_2 \xi_1 =  -  \partial_x \xi_1 , \\
 && \frac{1}{D_{23}} N_2 - \beta N_1 \xi_2 + \beta N_2 \xi_1 =  -  \partial_x \xi_2 ,
\end{eqnarray}
where $\alpha = \left(\frac{1}{D_{12}} - \frac{1}{D_{13}}\right)$, 
$\beta = \left(\frac{1}{D_{12}} - \frac{1}{D_{23}}\right)$.

We rewrite into:
\begin{eqnarray}
\label{ord_0}
&& \partial_t \xi_1 + \partial_x \cdot N_1 = 0 , \\
&& \partial_t \xi_2 + \partial_x \cdot N_2 = 0 , \\
&& \left( \begin{array}{c c}
\frac{1}{D_{13}} + \alpha \xi_2  & - \alpha \xi_1  \\
 - \beta \xi_2 &  \frac{1}{D_{23}}  + \beta \xi_1 
\end{array} \right)
 \left( \begin{array}{l}
N_1  \\
N_2 
\end{array} \right) = 
 \left( \begin{array}{l}
 -  \partial_x \xi_1 \\
 -  \partial_x \xi_2 
\end{array} \right)
\end{eqnarray}
and we have
\begin{eqnarray}
\label{part_0}
&& \partial_t \xi_1 + \partial_x \cdot N_1 = 0 , \\
\label{part_1}
&& \partial_t \xi_2 + \partial_x \cdot N_2 = 0 , \\
\label{part_2}
&& \left( \begin{array}{l}
N_1  \\
N_2 
\end{array} \right) = \frac{D_{13} D_{23}}{1 + \alpha D_{13} \xi_2 + \beta D_{23} \xi_1}
 \left( \begin{array}{c c}
\frac{1}{D_{23}} + \beta \xi_1  & \alpha \xi_1  \\
 \beta \xi_2 &  \frac{1}{D_{13}}  + \alpha \xi_2 
\end{array} \right)
 \left( \begin{array}{l}
 -  \partial_x \xi_1 \\
 -  \partial_x \xi_2 
\end{array} \right) .
\end{eqnarray}

The next step is to apply the semi-discretization of the 
partial differential operator $\frac{\partial}{\partial x}$.

We apply the first differential operator in equation (\ref{part_0}) and (\ref{part_1})
as an forward upwind scheme given as
\begin{eqnarray}
\frac{\partial}{\partial x} & = & D_+ =  \frac{1}{\Delta x}\cdot \left(\begin{array}{rrrrr}
 -1 & 0 & \ldots & ~ & 0 \\
  1 & -1 & 0 & \ldots & 0 \\
 \vdots & \ddots & \ddots & \ddots & \vdots \\
 0 & ~ & 1 & -1 & 0 \\
 0 & \ldots & 0 & 1 & -1
\end{array}\right)~\in~\R^{(J+1) \times (J+1)}
\end{eqnarray}
and  the second differential operator in equation (\ref{part_2})
as an backward upwind scheme given as
\begin{eqnarray}
\frac{\partial}{\partial x} & = & D_- =  \frac{1}{\Delta x}\cdot \left(\begin{array}{rrrrr}
 -1 & 1 & 0 & \ldots & 0 \\
 0 & -1 & 1 & 0 & \ldots \\
 \vdots & \ddots & \ddots & \ddots & \ddots \\
 0 & \ldots & 0  & -1 & 1 \\
 0 & ~ & \ldots & 0 & -1
\end{array}\right)~\in~\R^{(J+1) \times (J+1)} .
\end{eqnarray}

\subsection{Experiments with the Iterative scheme in time (Global Linearisation)}

In the first experiments, we test the first iterative scheme
(iterative scheme in time (Global Linearisation)).

\begin{figure}[ht]
\begin{center}  
\includegraphics[width=8.0cm,angle=-0]{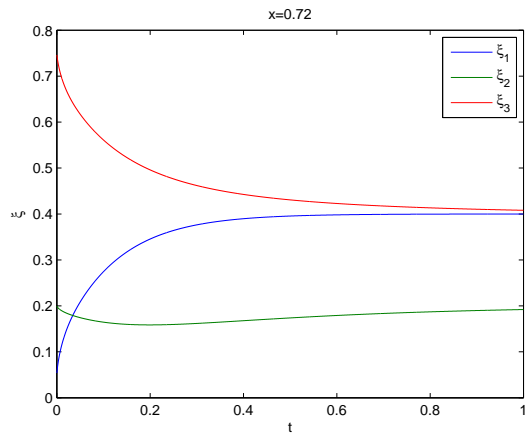}
\end{center}
\caption{\label{multi_1} The figures present the results of the
concentration $c_1$, $c_2$ and $c_3$.}
\end{figure}

We test the different schemes and obtain the results shown in 
Figure \ref{multi_1}.

\begin{figure}[ht]
\begin{center}  
\includegraphics[width=5.0cm,angle=-0]{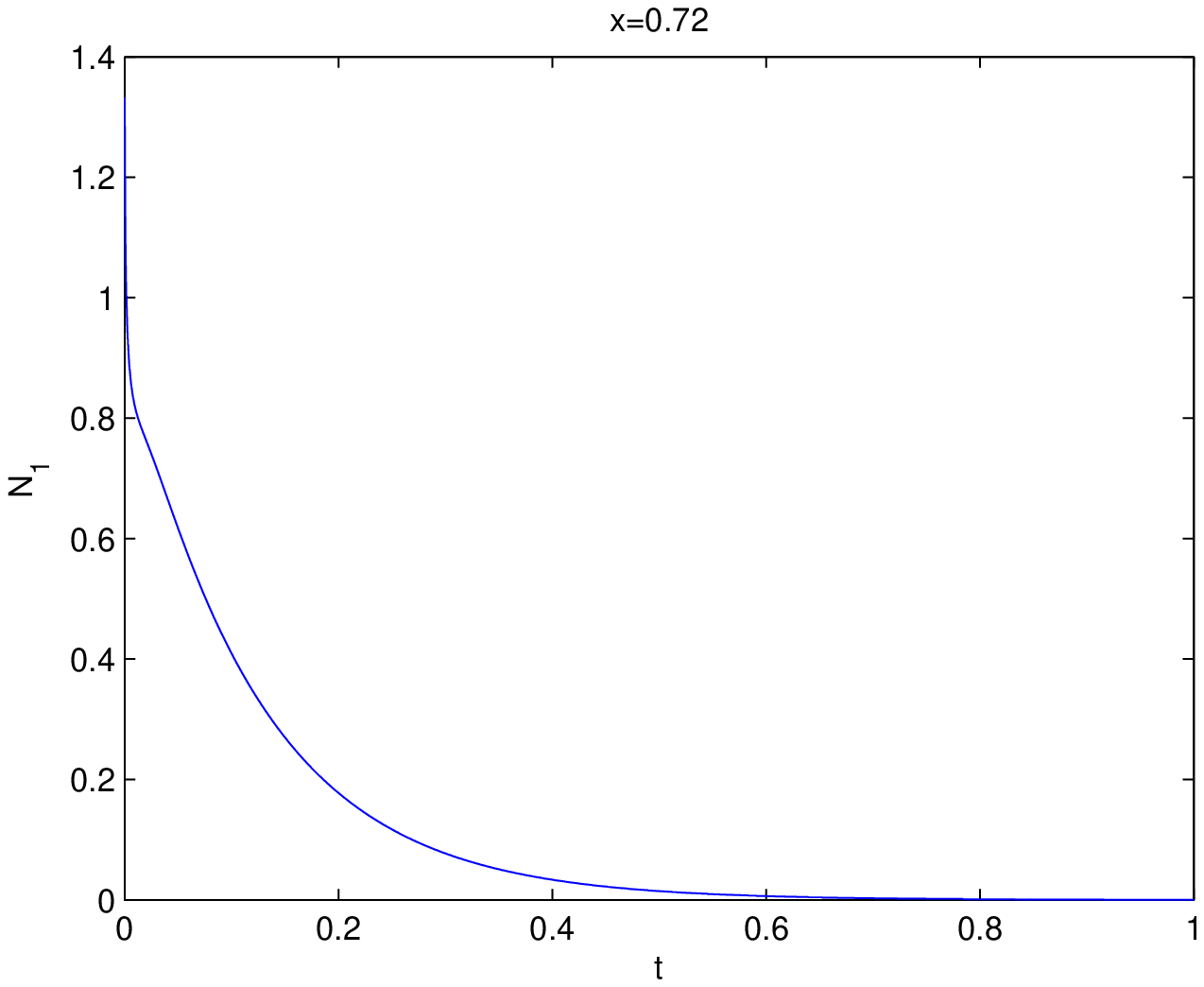}
\includegraphics[width=5.0cm,angle=-0]{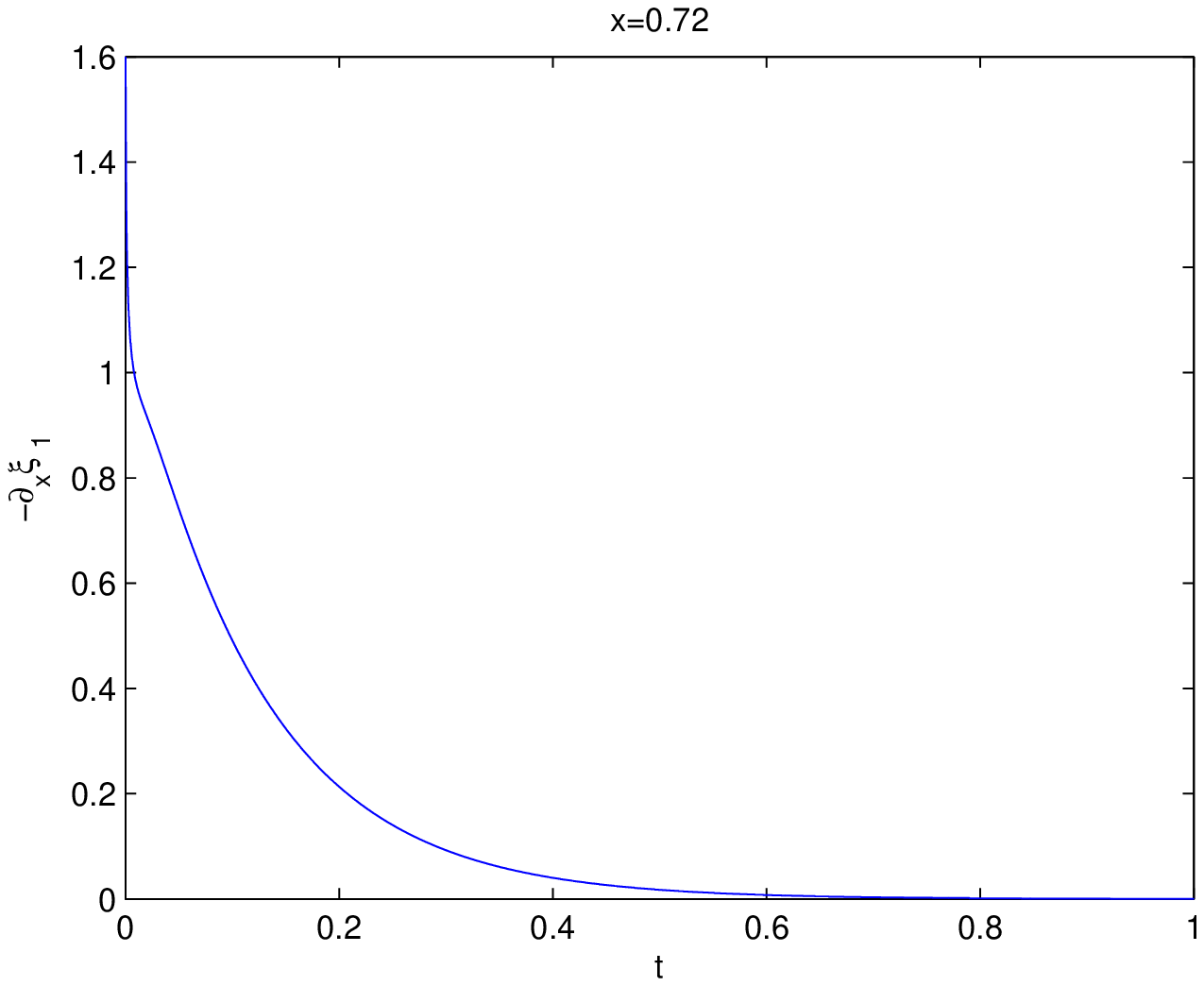} \\
\includegraphics[width=5.0cm,angle=-0]{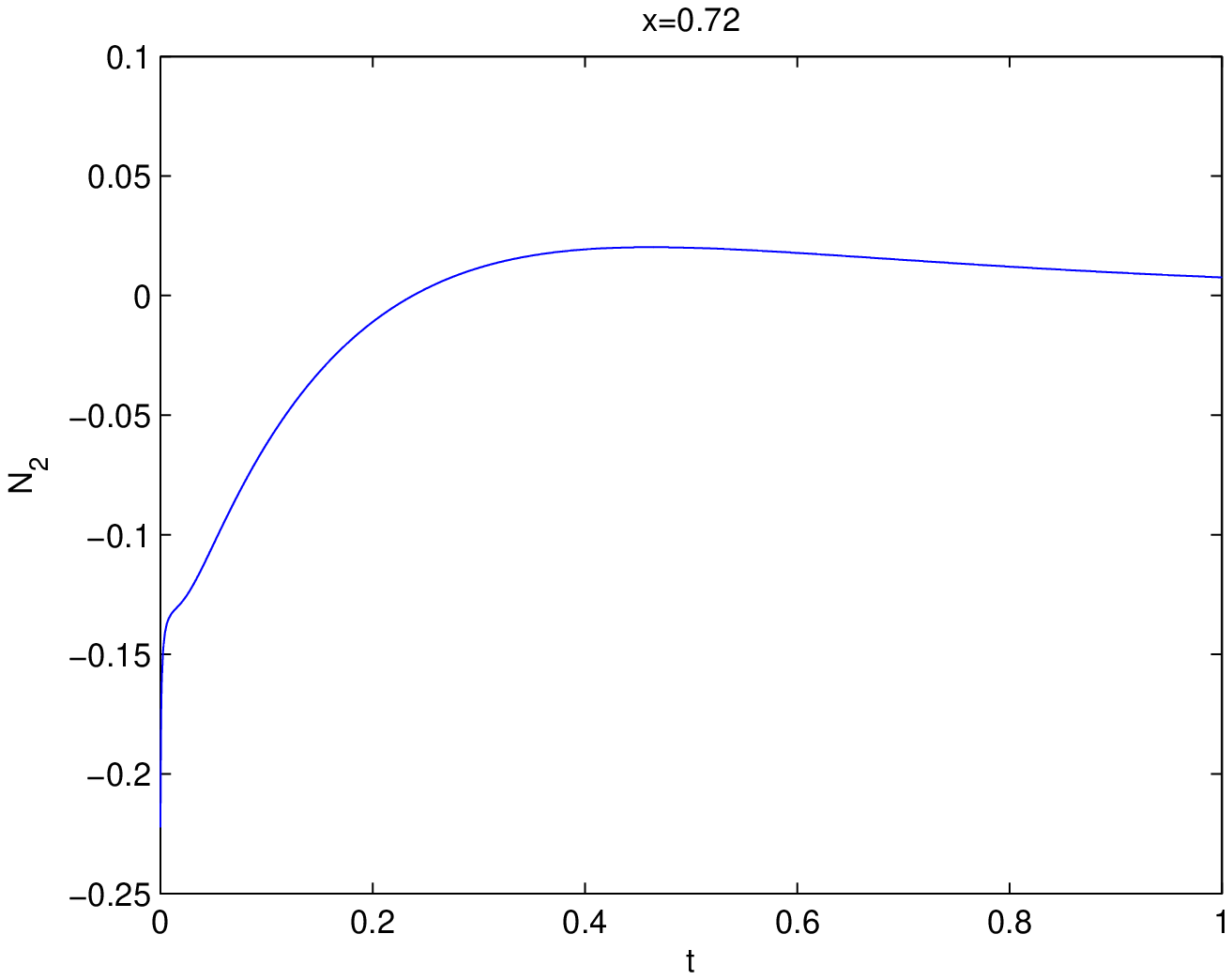}
\includegraphics[width=5.0cm,angle=-0]{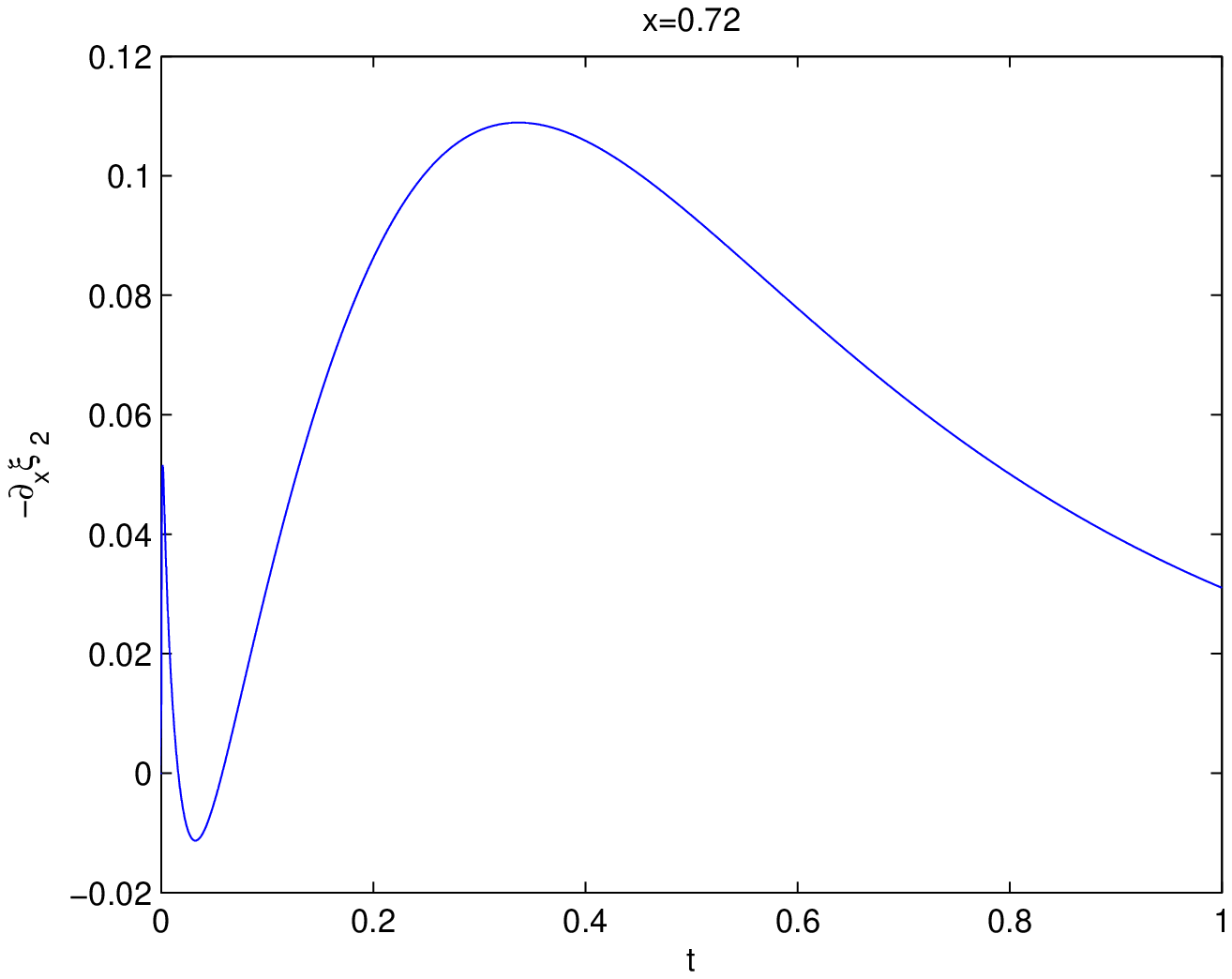}
\end{center}
\caption{\label{multi_2} The upper figures present the results of the
concentration $c_1$ and $- \partial_x \xi_1$. The lower figures presents
the results of $c_2$ and $- \partial_x \xi_2$.}
\end{figure}

The concentration and their fluxes are given in Figure \ref{multi_2}.

\begin{figure}[ht]
\begin{center}  
\includegraphics[width=5.0cm,angle=-0]{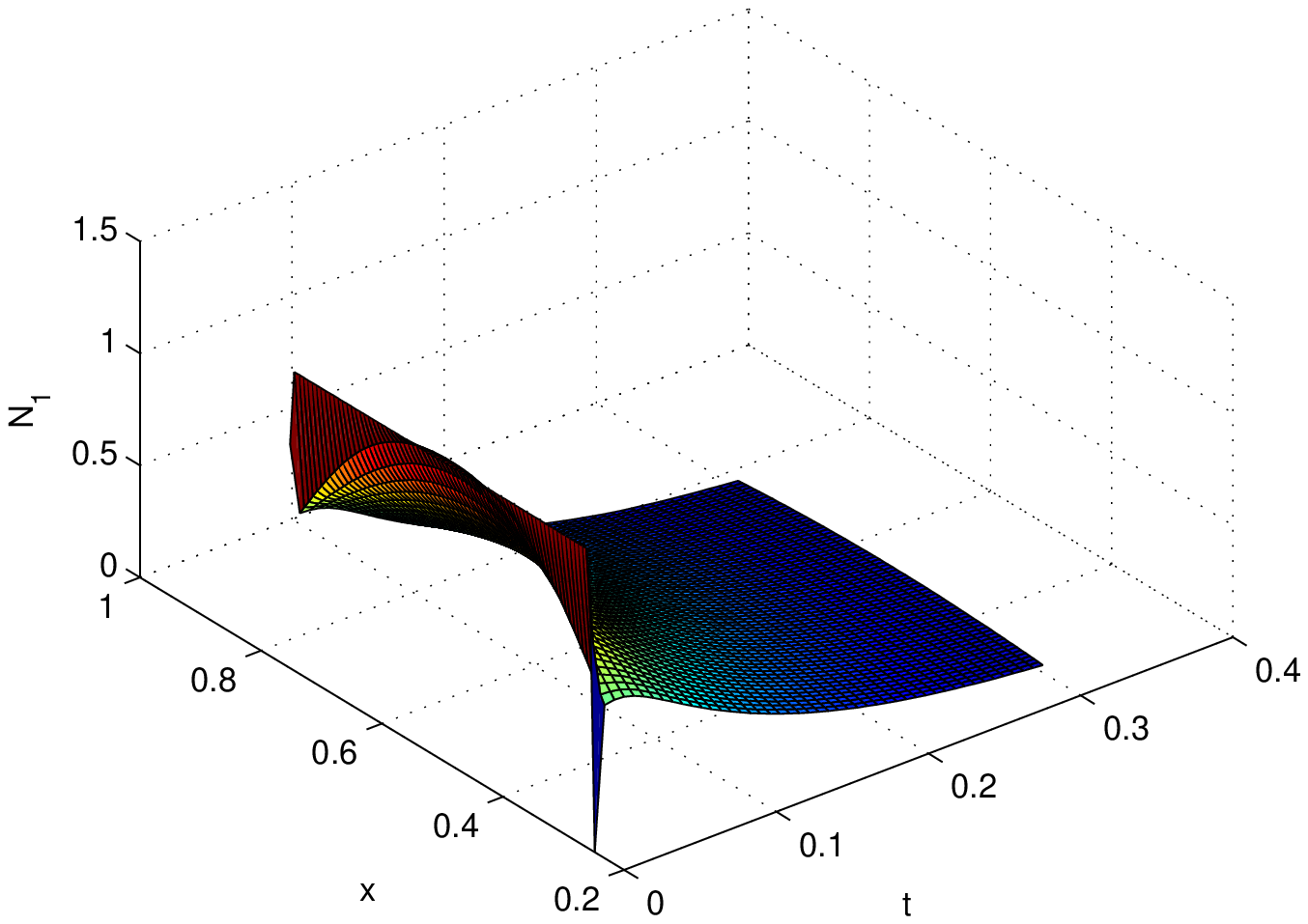}
\includegraphics[width=5.0cm,angle=-0]{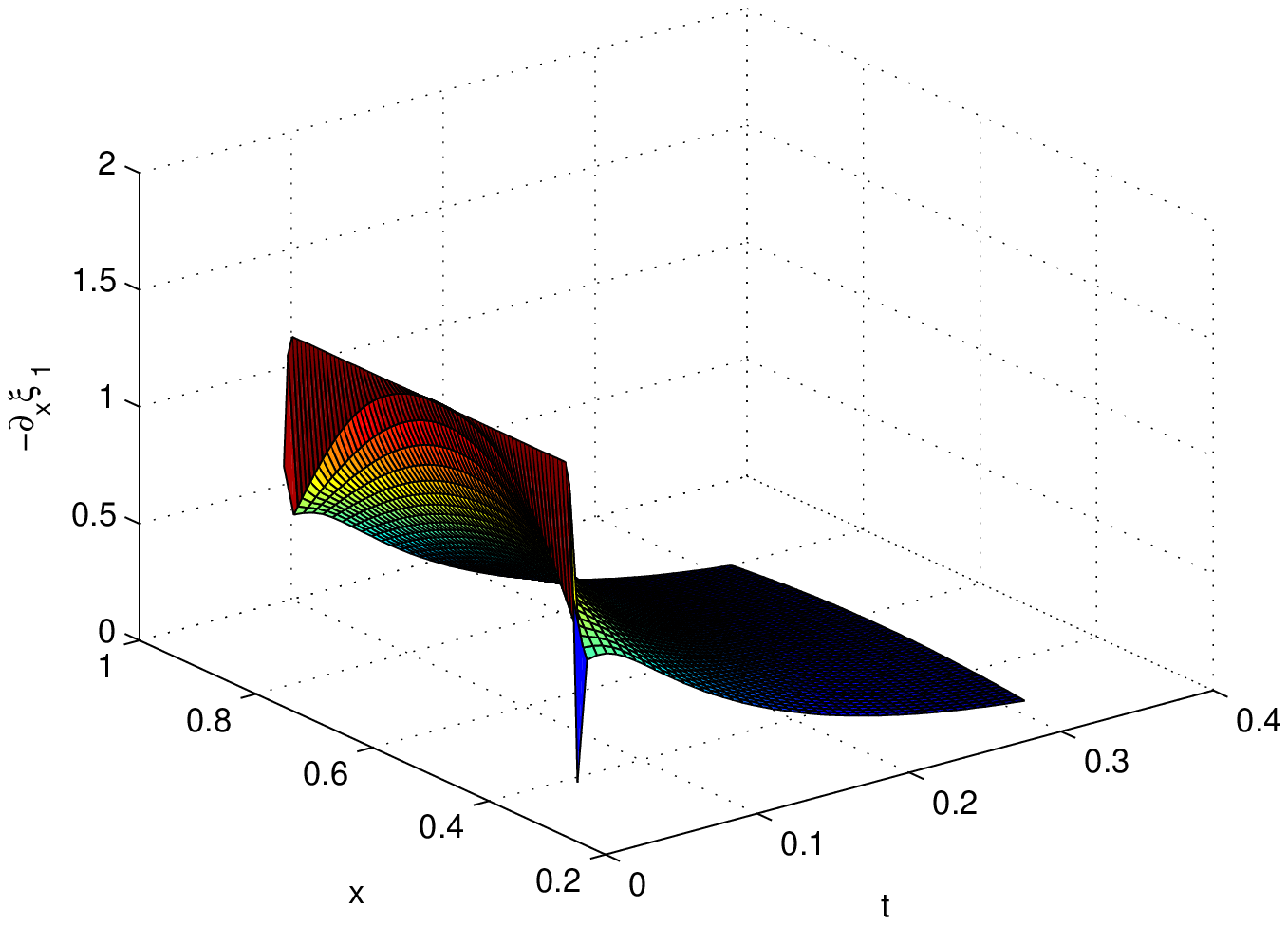} \\
\includegraphics[width=5.0cm,angle=-0]{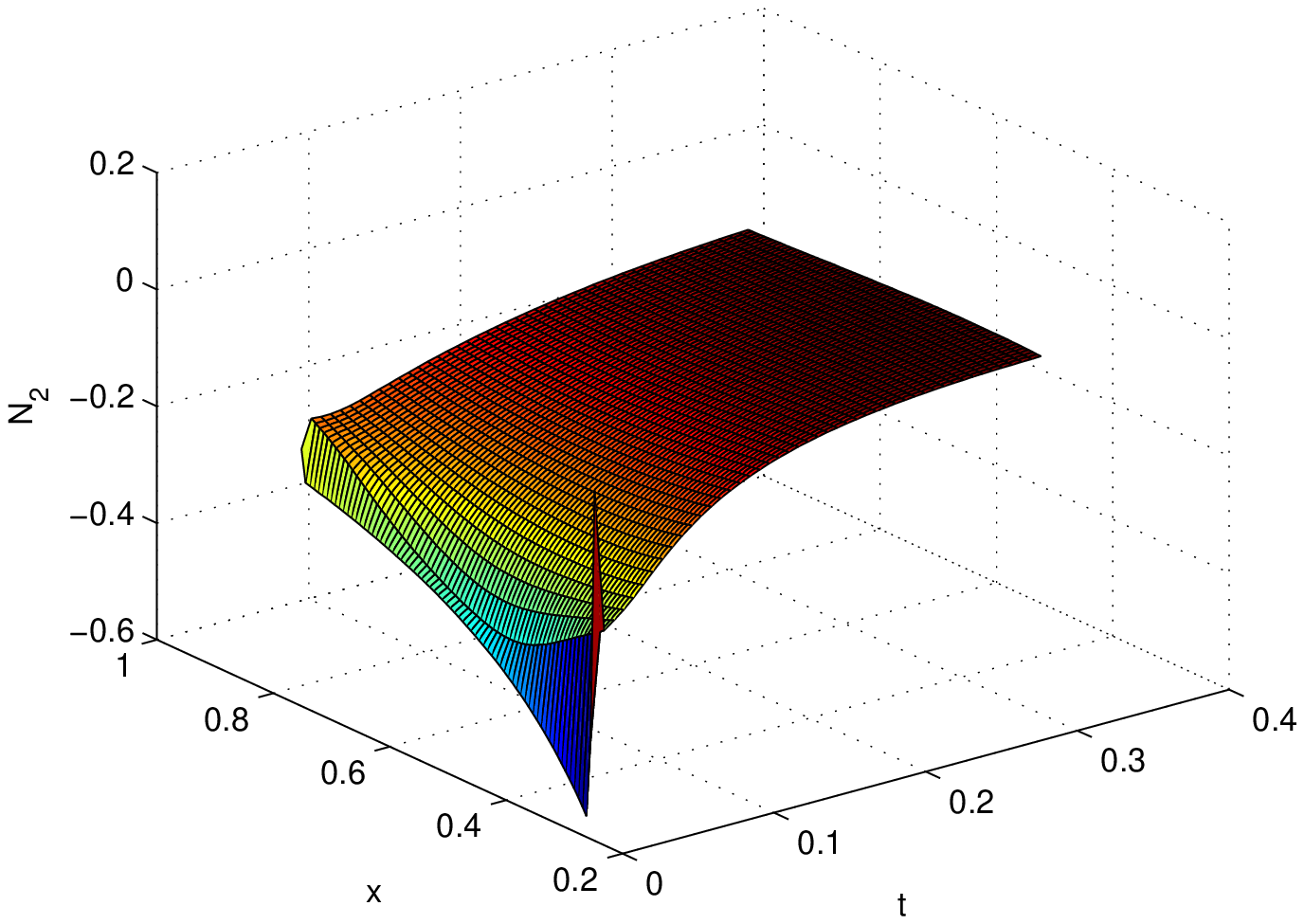}
\includegraphics[width=5.0cm,angle=-0]{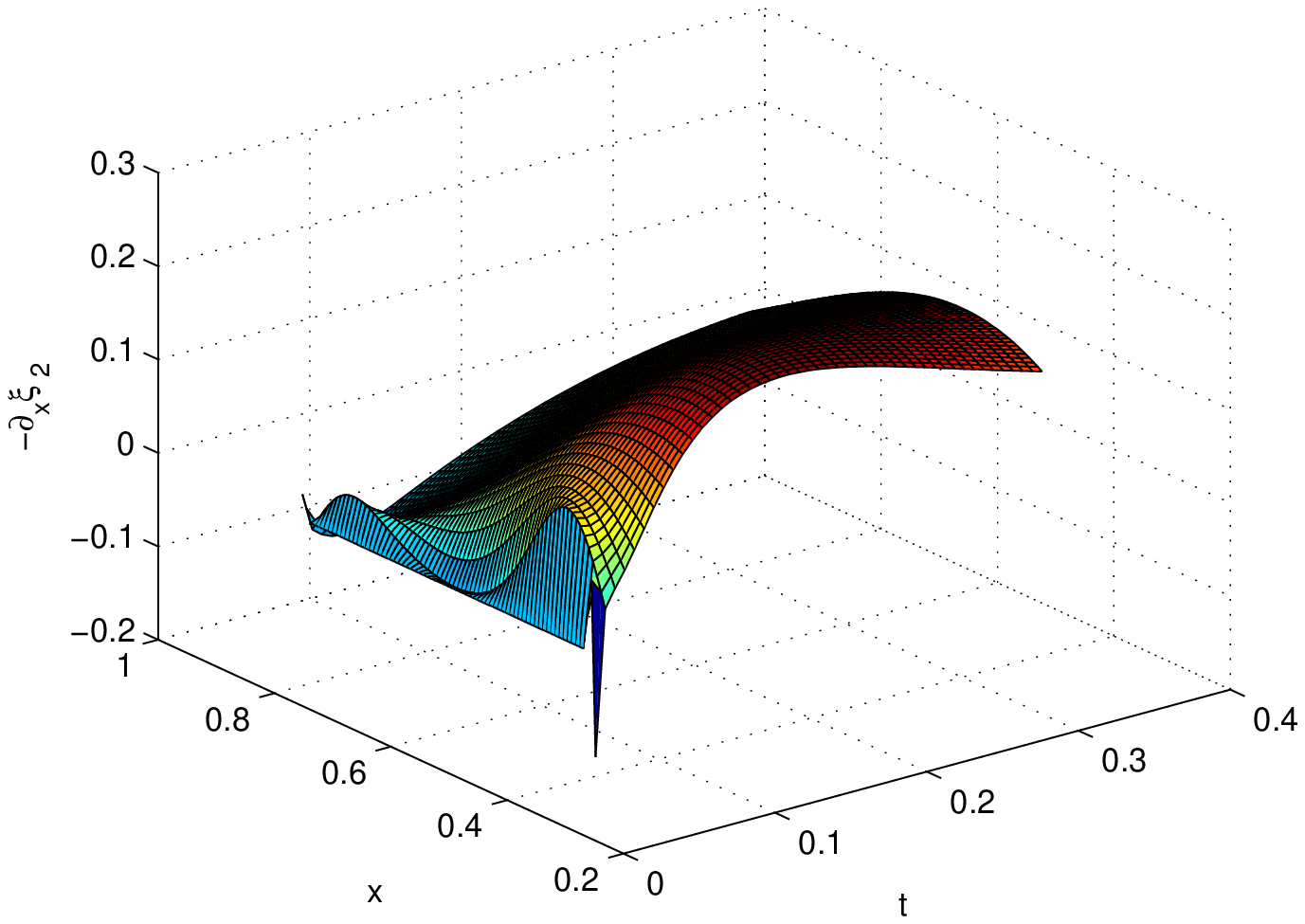}
\end{center}
\caption{\label{multi_3} The figures present the results of the 3d plots in time and space.  The upper figures present the results of the
concentration $c_1$ and $- \partial_x \xi_1$. The lower figures presents
the results of $c_2$ and $- \partial_x \xi_2$.}
\end{figure}

The full plots in time and space of the concentrations and their fluxes are given in Figure \ref{multi_3}.

\begin{figure}[ht]
\begin{center}  
\includegraphics[width=5.0cm,angle=-0]{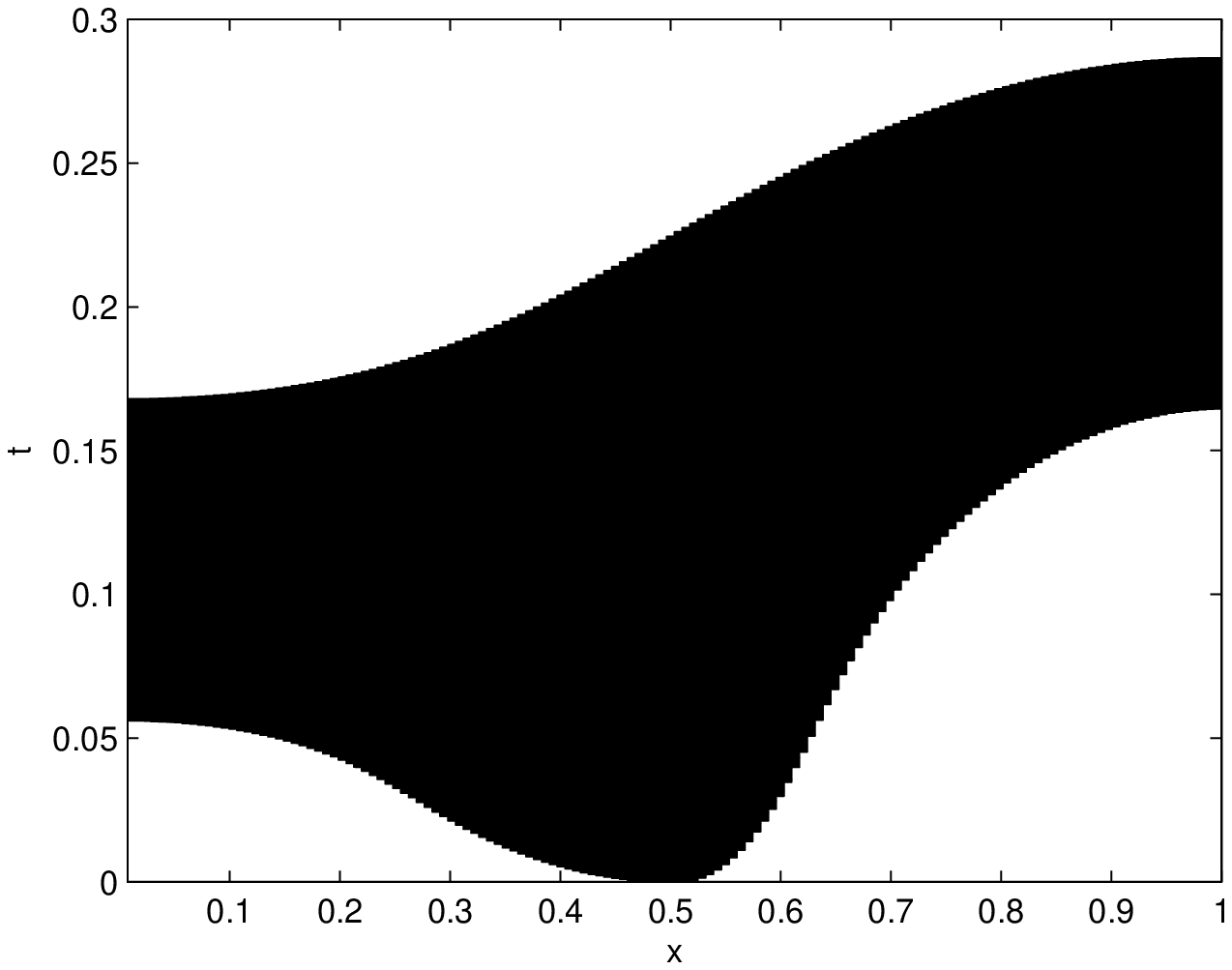}
\includegraphics[width=5.0cm,angle=-0]{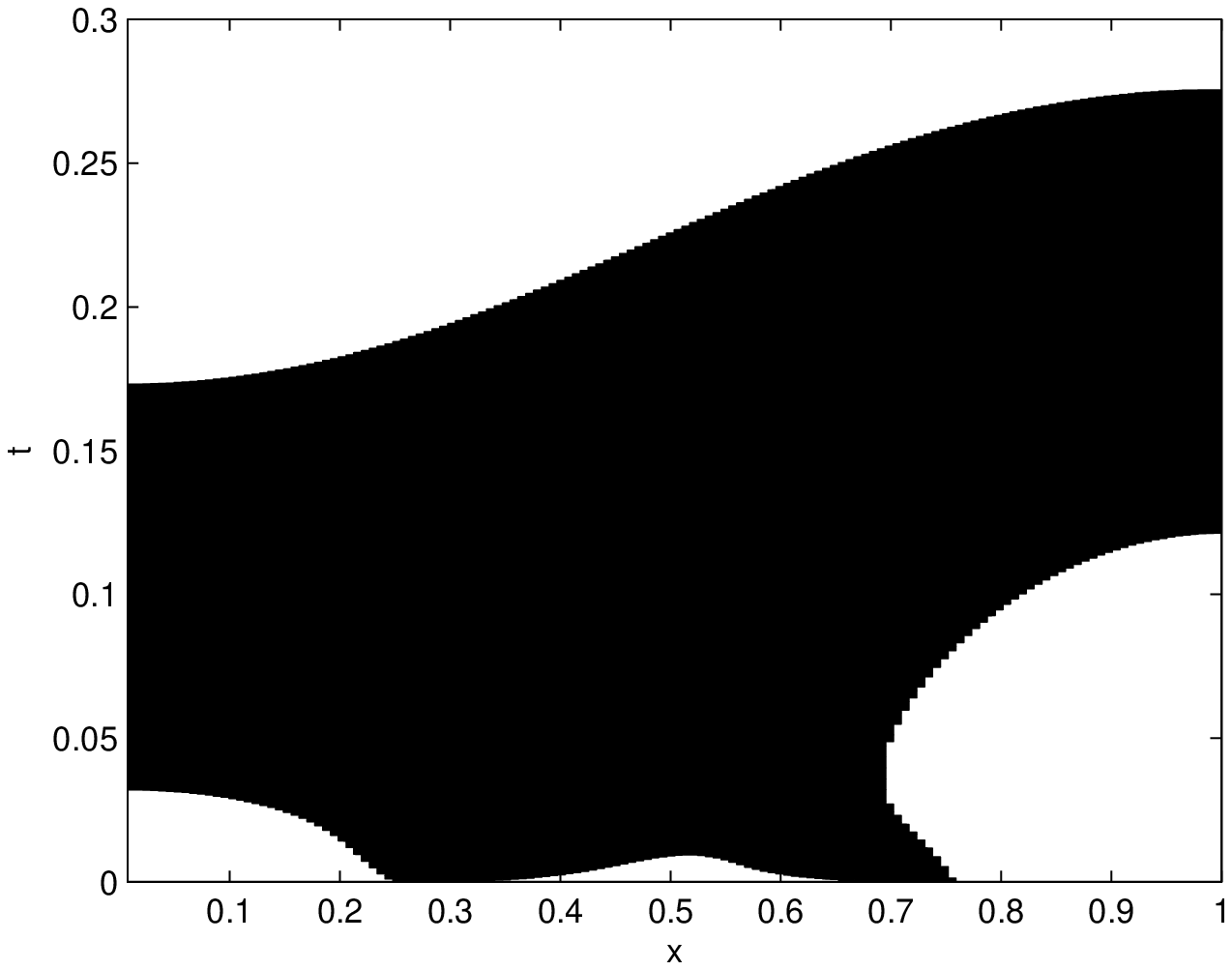} 
\end{center}
\caption{\label{multi_4} The figures present the asymptotic diffusion (left hand side) and uphill diffusion (right hand side) in the space-time region.}
\end{figure}

The space-time regions where $- N2 \partial_x \xi_2 \ge 0$ for the
uphill diffusion and asymptotic diffusion, given in Figure \ref{multi_4}.

\begin{remark}
The first method applies a global linearization based on the time-steps.
All effects are resolved, by the way, we have taken into account the
CFL condition. We achieve better results with finer time-steps, e.g.,
$\Delta t_{CFl} / 8$, such that the global linearisation, via the 
time-step is important.
\end{remark}

\subsection{Iterative Schemes in time (Local Linearisation)} 

In the next series of experiments, we apply the more refined
linearization scheme, means the iterative approximation in 
a single time-step.

We apply the numerical convergence of the schemes with the
reference solution of the explicit method by $\Xi_{ref} = (\xi_1, \xi_2)$
where the time-step is $\Delta t_{CFL} / 8$ for this refined solution the
error is only marginal.

Based on the reference solution, we deal with the following errors:
\begin{eqnarray}
\label{kap7_gleich3}
E_{L_{1, \Delta x}}(t) & = & \int_{\Omega} | \Xi_{method, J, \Delta x}(x, t) - \Xi_{ref}(x, t) | \; dx  \nonumber \\
              & = & \Delta x \sum_{i=1}^N | \Xi_{method, J, \Delta x}(x_i, t) - \Xi_{ref}(x_i, t) |  ,
\end{eqnarray}
where $method, J$ is the Richardson with $J$ iterative steps and \\
$\Delta t = \Delta t_{CFL}, \Delta t_{CFL}/2, \Delta t_{CFL}/4$.
Further $method, expl$ is the explicit method with $\Delta t = \Delta t_{CFL}, \Delta t_{CFL}/2, \Delta t_{CFL}/4$.

\begin{figure}[ht]
\begin{center}  
\includegraphics[width=5.0cm,angle=-0]{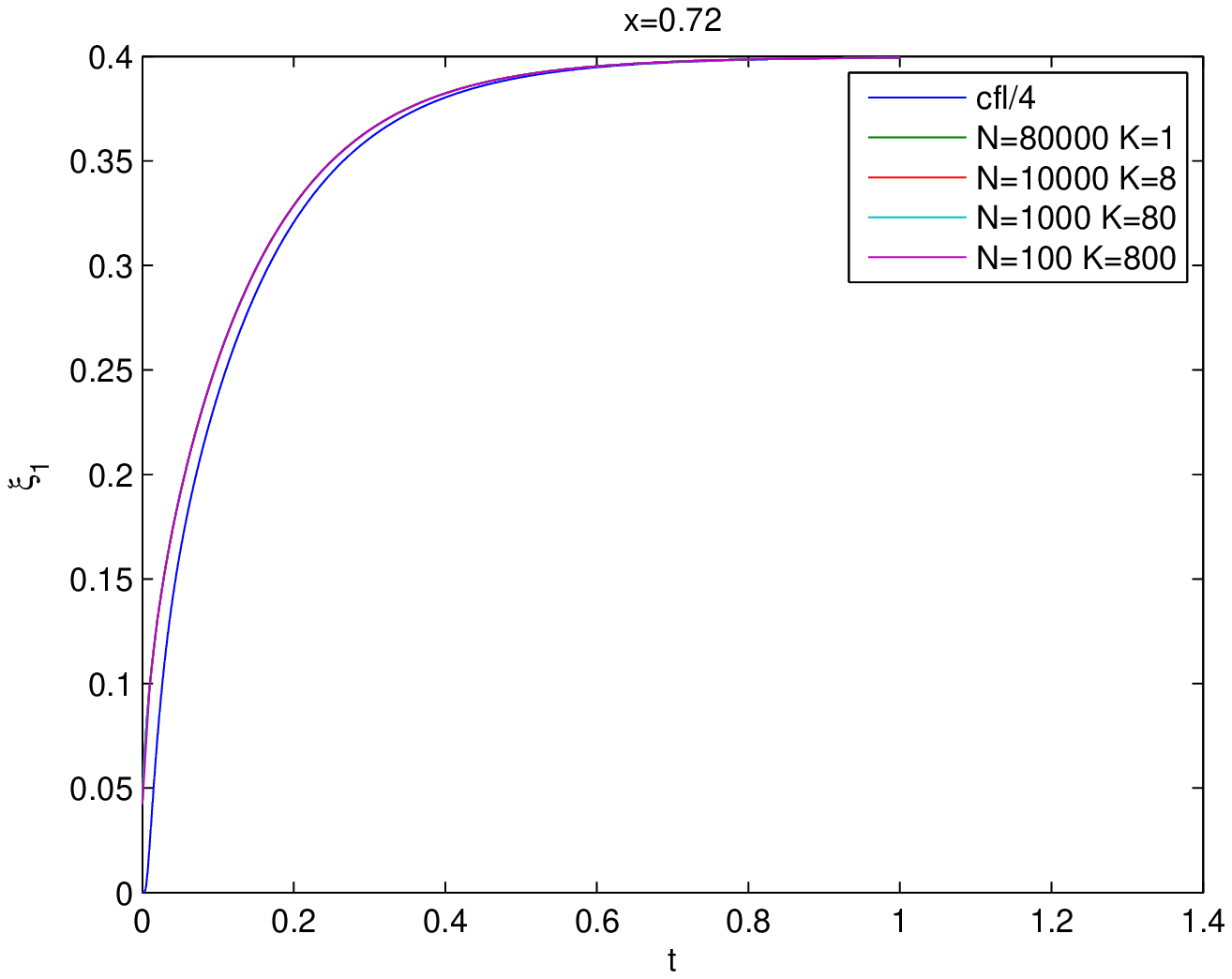}
\includegraphics[width=5.0cm,angle=-0]{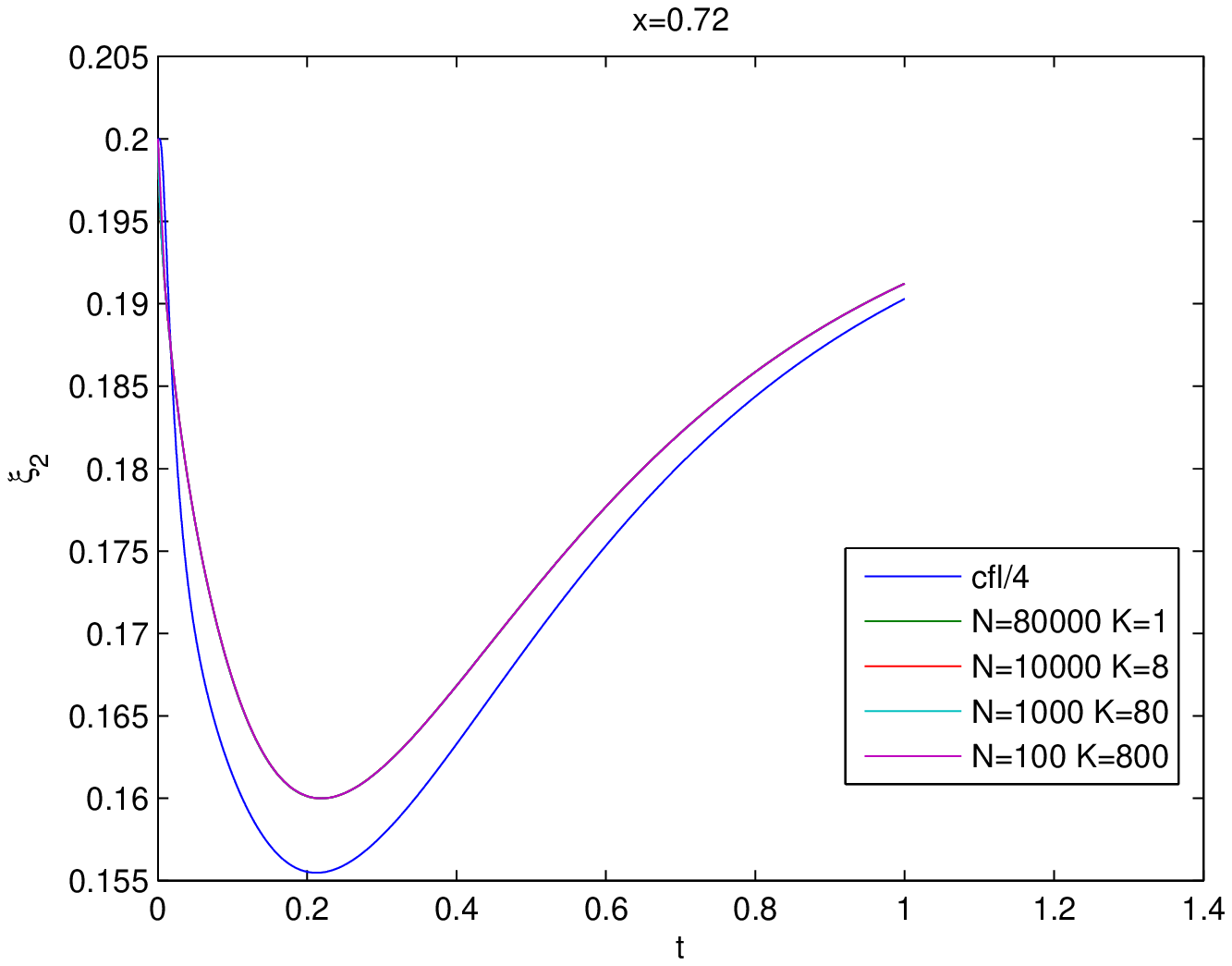} 
\end{center}
\caption{\label{multi_4_1}  The figures present the solutions of the
different time-step and iterative step of the Richardson-method}
\end{figure}

\begin{figure}[ht]
\begin{center}  
\includegraphics[width=5.0cm,angle=-0]{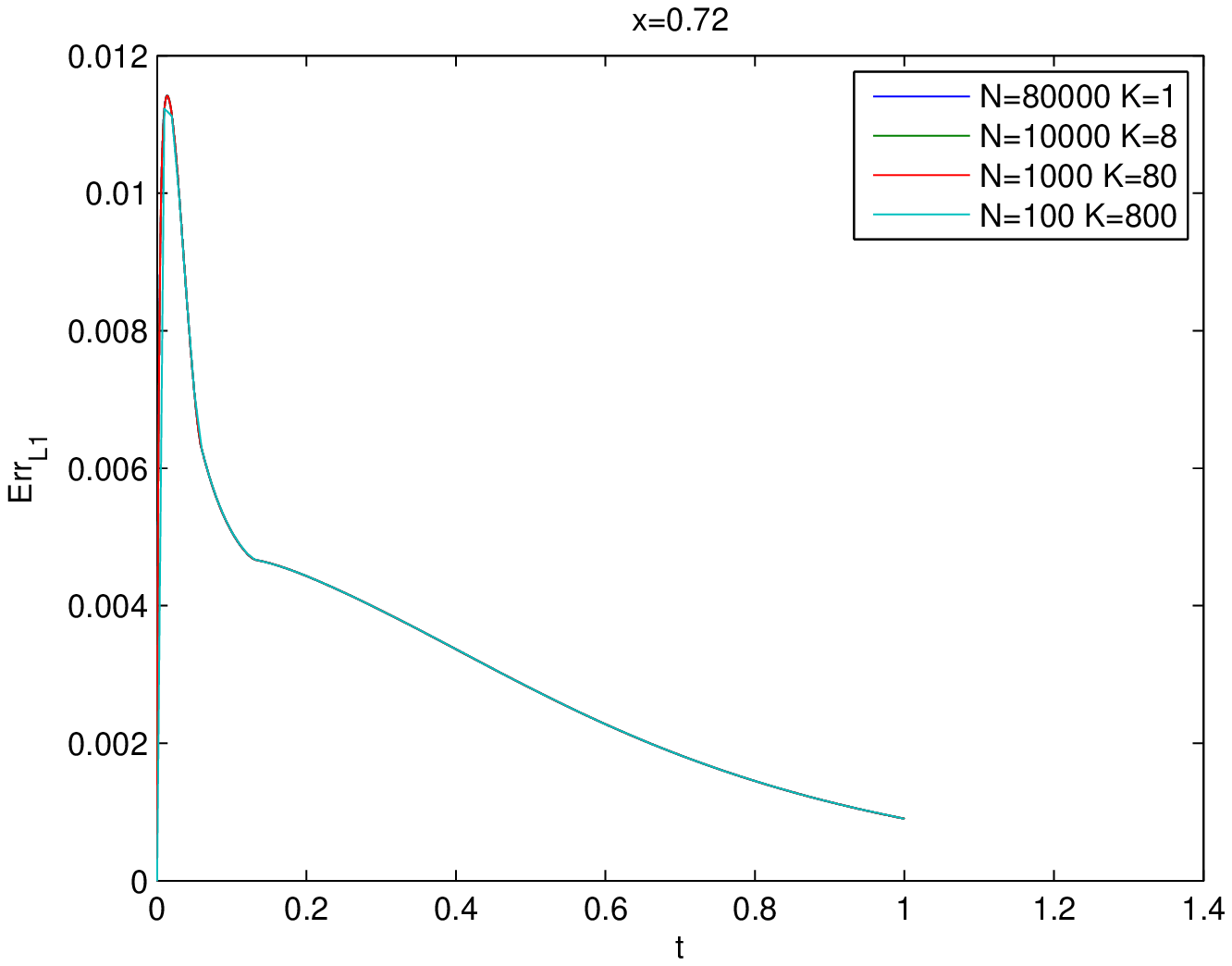}
\includegraphics[width=5.0cm,angle=-0]{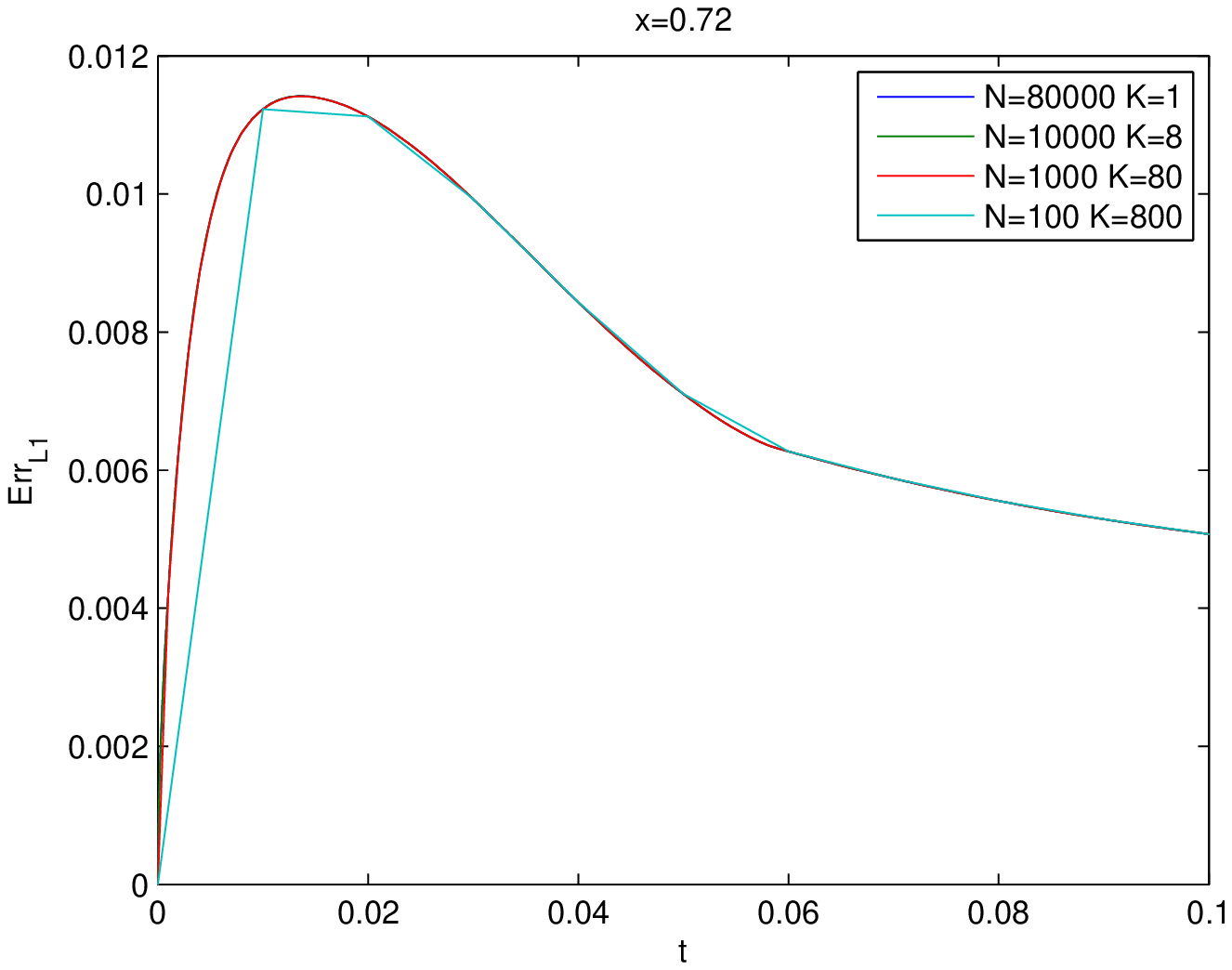} 
\end{center}
\caption{\label{multi_4_2} The figures present the errors of the
different time-step and iterative step solutions.}
\end{figure}

We apply the different versions of time-steps and iterative steps,
a reference solution is obtain with a fine time-step $\Delta t = \Delta t_{CFL}/4$. We see improvements in Figure \ref{multi_4_1} and the errors in Figure \ref{multi_4_2}.

\begin{remark}
Here, we see the benefit of large time-steps with $N=100$ and $K=800$,
means we have only $100$ time-steps and $800$ iterative steps, which are not expensive. Therefore we could gain the same results as with many small time-steps $N=80000$ and only one iterative step $K=1$. Such that the relaxation method benefits with the iterative cycles and we could enlarge the time-steps. 
\end{remark}

\begin{remark}
The second method applies a linear linearization based on the iterative
approaches in each single time-step.
We have the benefit of a relaxation in each local time-step,
such that we see a more accurate solution also with larger time-steps
than in the global linearization method.
\end{remark}

\section{Conclusions and Discussions }
\label{concl}

We present a fluid model based on Maxwell-Stefan diffusion equations.
The underlying problems for such a more delicate diffusion matrix is
discussed. Based on the nonlinear partial differential equations, we have to
apply linearisation approaches.
For first test-examples, we achieve more accurate results for a
so-called local linearized scheme. 
In future, we concentrate on the numerical convergence analysis and
generalize our results to real-life applications.

\bibliographystyle{plain}

\end{document}